\documentclass[reqno]{amsart}

\usepackage[latin1]{inputenc}
\usepackage[T1]{fontenc}                  
\usepackage[english]{babel}
\usepackage{amsmath,amssymb,amsbsy,amsthm,amsfonts}
\usepackage{color}
\usepackage{epsfig}
\usepackage{pgf}

\newtheorem{defin}{Definition}[section]

\newtheorem{remark}[defin]{Remark}

\newtheorem{theorem}[defin]{Theorem}

\usepackage[colorlinks=true]{hyperref}
\usepackage{theoremref}
\newcommand{\bluee}{\textcolor{black}}

\usepackage{accents}

\usepackage{wasysym}

\usepackage{enumerate}
\usepackage{xcolor}

\begin{document}
\title[Nonlinear estimates for reaction diffusion equations] 
      {Nonlinear estimates for traveling wave solutions of reaction diffusion equations}  
\author[L.-C. Hung and X. Liao]{Li-Chang Hung$^{\ast}$ and Xian Liao$^{\natural}$}

\email{lichang.hung@gmail.com; xian.liao@kit.edu}
\thanks{$^{\ast}$Department of Mathematics, National Taiwan University, Taipei, Taiwan}
\thanks{$^{\natural}$Institute of Analysis, Karlsruhe Institute of Technology, Karlsruhe, Germany}


\begin{abstract}
In this paper we will establish nonlinear a priori lower and upper bounds for the solutions to a large class of equations which arise from the study of traveling wave solutions of reaction-diffusion equations, and we will apply our nonlinear bounds to the Lotka-Volterra system of two competing species as examples.
The idea used in a series of papers \cite{NBMP-Discrete,JDE-16,CPAA-16,DCDS-B-18,NBMP-n-species,DCDS-A-17} for the establishment of the linear N-barrier maximum principle  will also be used in the proof. 

\end{abstract}

\maketitle

\section{Introduction}\label{sec: intro}
The present paper is devoted to \textit{nonlinear} a priori upper and lower bounds for the solutions $u_i=u_i(x): \mathbb{R}\mapsto[0,\infty)$, $i=1,\cdots, n$ to the following boundary value problem of $n$  equations 
\begin{equation}\label{eqn}
\begin{cases}
\vspace{3mm}
d_i\,(u_i)_{xx}+\theta\,(u_i)_{x}+u_i^{l_i}\,f_i(u_1,u_2,\cdots,u_n)=0, \ \ x\in\mathbb{R}, \ \ i=1,2,\cdots,n, \\
(u_1,u_2,\cdots,u_n)(-\infty)=\textbf{e}_{-},\quad (u_1,u_2,\cdots,u_n)(\infty)=\textbf{e}_{+}.
\end{cases}
\end{equation} 
In the above,   $d_i$, $l_i>0$, $\theta\in\mathbb{R}$ are parameters, $f_i\in C^0([0,\infty)^n)$ are given functions   and the boundary values $\textbf{e}_{-}, \textbf{e}_{+}$ take value in the following constant equilibria  set
\begin{equation}\label{eqn: e- and e+}
\Big\{ (u_1,\cdots,u_n) \;\Big|\; u_i^{l_i}\,f_i(u_1,\cdots,u_n)=0,\quad u_i\geq 0,\quad \forall i=1,\cdots,n\Big\}.
\end{equation}

Equations \eqref{eqn} arise from the study of traveling waves solutions of reaction-diffusion equations (see \cite{Murray93Mbiology,Volpert94}).
A series of papers  \cite{NBMP-Discrete,JDE-16,CPAA-16,DCDS-B-18,NBMP-n-species,DCDS-A-17} by Hung \textit{et al.} have been contributed to the \textit{linear} (N-barrier) maximum principle for  the $n$ equations \eqref{eqn}, and in particular the   lower and upper bounds for any linear combination of the solutions
$$\displaystyle\sum_{i=1}^{n}\alpha_i\,u_i(x),\quad \forall (\alpha_1,\cdots,\alpha_n)\in (\mathbb{R}^+)^n$$ 
have been established in terms of \bluee{the parameters $d_i,l_i,\theta$ in \eqref{eqn}}. 






Here we aim to derive \textit{nonlinear} estimates for the polynomials of the solutions:
 $$\displaystyle\prod_{i=1}^{n} (u_i(x)+k_i)^{\alpha_i},\quad \forall(\alpha_1,\cdots,\alpha_n)\in(\mathbb{R}^+)^n$$
  for some $k_i\ge0$, which is related to the diversity indices of the species in ecology: $D^q=(\sum_{i=1}^n (u_i)^q)^{1/(1-q)}$, \bluee{$q\in [1,\infty)$}. 
Observe that when either $\text{\bf e}_{+}=(0,\cdots,0)$ or $\text{\bf e}_{-}=(0,\cdots,0)$, the trivial lower bound of $\displaystyle\prod_{i=1}^{n} (u_i(x)+k_i)^{\alpha_i}$ is $\prod_{i=1}^n k_i^{\alpha_i}$. 
For $k_i>0$ the following lower bound for the upper solutions of \eqref{eqn} holds.
\begin{theorem} [\textbf{Lower bound}]\thlabel{prop: lower bed}
Suppose that $(u_i(x))_{i=1}^n\in (C^2(\mathbb{R}))^n$ with $u_i(x)\ge0$, $\forall i=1,\cdots,n$ is an upper solution of \eqref{eqn}: 
\begin{equation}\label{eqn:upper} 
\begin{cases}
\vspace{3mm}
d_i(u_i)_{xx}+\theta(u_i)_{x}+u_i^{l_i}f_i(u_1,u_2,\cdots,u_n)\le0, \quad x\in\mathbb{R}, \quad i=1,2,\cdots,n, \\
(u_1,u_2,\cdots,u_n)(-\infty)=\textbf{e}_{-},\quad (u_1,u_2,\cdots,u_n)(\infty)=\textbf{e}_{+},
\end{cases}
\end{equation}
and that there exist $(\underaccent\bar{u}_i)_{i=1}^n\in(\mathbb{R}^+)^n$ such that 
\begin{equation}\label{H1}
\begin{array}{l} 
f_i(u_1,\cdots,u_n)\ge 0, \hbox{ for all }\\
 (u_1, \cdots,u_n)\in \underaccent\bar{\mathcal{R}}:=\big\{ (u_i)_{i=1}^n\in([0,\infty))^n\; |\; \sum_{i=1}^{n}\frac{\displaystyle u_i}{\displaystyle  \underaccent\bar{u}_i}\le 1\big\}.
\end{array}
\end{equation} 
Then we have for any $(k_i)_{i=1}^n\in (\mathbb{R}^+)^n$ and $(\alpha_i)_{i=1}^n\in (\mathbb{R}^+)^n$, 
\begin{equation}\label{ineq:lower}
\prod_{i=1}^{n} (u_i(x)+k_i)^{\bluee{d_i}\alpha_i}\ge e^{\lambda_1}, \quad x\in\mathbb{R},
\end{equation}
where
\begin{subequations}\label{eqn: lambda1 eta lambda2 lower bound}
\begin{alignat}{4}
\label{eqn: lambda1 lower bound}
\lambda_1&=\min_{1\le j \le n} \Bigl(\eta\,d_j+\sum_{i=1, i\neq j}^{n}\alpha_i(d_i-d_j)\ln k_i\Bigr),\\
\label{eqn: eta lower bound}
\eta&=\min_{1\le j \le n} \frac{1}{d_j}\Bigl(\lambda_2-\sum_{i=1, i\neq j}^{n}\alpha_i(d_i-d_j)\ln k_i\Bigr),\\
\label{eqn: lambda2 lower bound}
\lambda_2&=\min_{1\le j \le n} \Bigl(\alpha_j d_j \ln (\underaccent\bar{u}_j+k_j)+\sum_{i=1, i\neq j}^{n}\alpha_i d_i\ln k_i\Bigr).
\end{alignat}
\label{straincomponent}
\end{subequations}

\end{theorem}

\begin{remark}[Equal diffusion]\thlabel{rem: lower bound equal diffusion}
When $d_i=d$ for all $i=1,2,\cdots,n$, then
\begin{equation*}
\lambda_1=\min_{1\le j \le n} \Bigl(\alpha_j \ln (\underaccent\bar{u}_j+k_j)+\sum_{i=1, i\neq j}^{n}\alpha_i\ln k_i\Bigr)d=\lambda_2=d\eta,
\end{equation*}  
and
 the lower bound \eqref{ineq:lower} becomes
\begin{equation*}\label{eqn: lower bound of p equal diffusion}
\prod_{i=1}^{n} (u_i(x)+k_i)^{\alpha_i}
\ge \min_{1\leq j\leq n}\Bigl( (\underaccent\bar{u}_j+k_j)^{\alpha_j}\prod_{i\neq j}k_i^{\alpha_i}\Bigr), \quad x\in\mathbb{R}.
\end{equation*} 
If furthermore $\alpha_i=\alpha$, $\forall i=1,\cdots,n$, then  the inequality of arithmetic and geometric averages yields
\begin{align*}
\sum_{i=1}^n (u_i+k_i)^{\alpha}\geq n\Bigl(\prod_{i=1}^n (u_i+k_i)^{\alpha}\Bigr)^{\frac1n}
\geq n\min_{1\leq j\leq n}\Bigl( (\underaccent\bar{u}_j+k_j)^{\alpha}\prod_{i\neq j}k_i^{\alpha}\Bigr)^{\frac1n}.
\end{align*}
\end{remark}

On the other hand, we can find an upper bound of $\displaystyle\prod_{i=1}^{n} (u_i(x))^{\alpha_i}$ for the lower solutions of \eqref{eqn}. 
\begin{theorem} [\textbf{Upper bound}]\thlabel{prop: upper bed}
Suppose that $(u_i(x))_{i=1}^n\in (C^2(\mathbb{R}))^n$ with $u_i(x)\ge0$ $\forall i=1,\cdots,n$ is a lower solution of \eqref{eqn}:
\begin{equation}\label{eqn:lower} 
\begin{cases}
\vspace{3mm}
d_i(u_i)_{xx}+\theta(u_i)_{x}+u_i^{l_i}f_i(u_1,u_2,\cdots,u_n)\ge0, \quad x\in\mathbb{R}, \quad i=1,2,\cdots,n, \\
(u_1,u_2,\cdots,u_n)(-\infty)=\textbf{e}_{-},\quad (u_1,u_2,\cdots,u_n)(\infty)=\textbf{e}_{+},
\end{cases}
\end{equation}
and there exist $\bar{u}_i>0, i=1,\cdots,n,$ such that
\begin{equation}\label{H2}
\begin{array}{l} 
f_i(u_1,\cdots,u_n)\le 0,\hbox{ for all }\\
 (u_1,\cdots,u_n)\in \bar{\mathcal{R}}:=\big\{ (u_i)_{i=1}^n\in([0,\infty))^n\; |\; \sum_{i=1}^{n}\frac{\displaystyle u_i}{\displaystyle\bar{u}_i}\ge 1 \big\}.
\end{array}
\end{equation} 

Then we have for any $m_i\geq 1$ and $\alpha_i>0$ $(i=1,2,\cdots,n)$
\begin{equation}\label{bound:upper}
\sum_{i=1}^{n} \alpha_i (u_i(x))^{m_i}
\le
\left(\max_{1\le i \le n} \alpha_i (\bar{u}_i)^{m_i}\right) 
\frac{\displaystyle\max_{1\le i \le n} d_i}{\displaystyle\min_{1\le i \le n} d_i}, \quad x\in\mathbb{R},
\end{equation}
and hence
\begin{equation}\label{eqn: upper bound of p general case}
\prod_{i=1}^{n} (u_i(x))^{m_i/n}\leq 
\frac{\displaystyle\max_{1\le i \le n}\alpha_i\,\bar{u}_i^{m_i}}{n \left(\displaystyle\prod_{i=1}^{n} \alpha_i\right)^{1/n}}
\frac{\displaystyle\max_{1\le i \le n} d_i}{\displaystyle\min_{1\le i \le n} d_i}, \quad x\in\mathbb{R}. 
\end{equation}
In particular, when $\alpha_i=\alpha$ for all $i=1,\cdots,n$, \eqref{eqn: upper bound of p general case} becomes
\begin{equation}\label{eqn: upper bound of p}
\prod_{i=1}^{n} (u_i(x))^{m_i/n}\leq 
\frac{\displaystyle\max_{1\le i \le n}\bar{u}_i^{m_i}}{n}
\frac{\displaystyle\max_{1\le i \le n} d_i}{\displaystyle\min_{1\le i \le n} d_i}, \quad x\in\mathbb{R}.
\end{equation}

\end{theorem}

In order to prove \thref{prop: lower bed}, we will first rewrite the system \eqref{eqn:upper} into the system for the new unknowns $(U_i)_{i=1}^n:=(\ln(u_i+k_i))_{i=1}^n$. Then we will follow the ideas in  \cite{NBMP-Discrete,JDE-16,CPAA-16,DCDS-B-18,NBMP-n-species,DCDS-A-17} to establish the lower bound for the linear combination of $(U_i)_{i=1}^n$, which implies the nonlinear lower bound \eqref{ineq:lower} correspondingly. Similarly, we will consider the new unknowns $(U_i)_{i=1}^n:=(u_i^{m_i})_{i=1}^n$ to establish the upper bound \eqref{bound:upper}. The proofs will be found in Section~\ref{Sec: Proofs}.


As an example to illustrate our main result, we use the Lotka-Volterra system of two competing species to conclude with Section~\ref{sec: intro}. This example provides an intuitive idea of the construction of the N-barrier in multi-species cases.

To illustrate \thref{prop: upper bed} for the case $n=2$, we use the Lotka-Volterra system of two competing species coupled with Dirichlet boundary conditions: 
\begin{equation}\label{eqn: 2 species LV rem}
\begin{cases}
\vspace{3mm} 
d_1 u_{xx}+\theta\,u_x+u\,(1-u-a_1\,v)=0,\ \ &x\in\mathbb{R},\\
\vspace{3mm} 
d_2 v_{xx}+\theta\,v_x+\kappa\,v\,(1-a_2\,u-v)=0,\ \ &x\in\mathbb{R},\\
(u,v)(-\infty)=\textbf{e}_i,\ \ (u,v)(+\infty)=\textbf{e}_j,
\end{cases}
\end{equation}
where $a_1$, $a_2$, $\kappa>0$ are constants. In \eqref{eqn: 2 species LV rem}, the constant equilibria are $\textbf{e}_1=(0,0)$, $\textbf{e}_2=(1,0)$, $\textbf{e}_3=(0,1)$ and $\textbf{e}_4=(u^{\ast},v^{\ast})$, where $(u^{\ast},v^{\ast})=\left(\displaystyle\frac{1-a_1}{1-a_1\,a_2},\displaystyle\frac{1-a_2}{1-a_1\,a_2}\right)$ is the intersection of the two straight lines $1-u-a_{1}\,v=0$ and $1-a_{2}\,u-v=0$ whenever it exists. We call the solution $(u(x),v(x))$ of \eqref{eqn: 2 species LV rem} an $(\textbf{e}_i,\textbf{e}_j)$-wave.

Tang and Fife (\cite{tang1980pfc}), and Ahmad and Lazer (\cite{Ahmad&Lazer91}) established the existence of the $(\textbf{e}_1,\textbf{e}_4)$-waves. 
Kan-on (\cite{Kan-on95,Kan-on97Fisher-Monostable}), Fei and Carr (\cite{Fei&Carr03}), Leung, Hou and Li (\cite{Leung08}), and Leung and Feng (\cite{Leung&Feng09}) proved the existence of $(\textbf{e}_2,\textbf{e}_3)$-waves using different approaches. $(\textbf{e}_2,\textbf{e}_4)$-waves were studied for instance, by Kanel and Zhou (\cite{Kanel&Zhou96}), Kanel (\cite{Kanel06}), and Hou and Leung (\cite{Hou&Leung08}).

For the above-mentioned $(\textbf{e}_1,\textbf{e}_4)$-waves, $(\textbf{e}_2,\textbf{e}_3)$-waves, and $(\textbf{e}_2,\textbf{e}_4)$-waves, we show a lower and an upper bounds of $u(x)v(x)$ by \thref{prop: lower bed} and \thref{prop: upper bed} respectively.
 To this end, let 
\begin{align*}
&\underline{u}=\min\left(1,\frac{1}{a_2}\right),\quad \bar{u} = \max \left(1,\frac{1}{a_2}\right),\\ 
&\underline{v}=\min\left(1,\frac{1}{a_1}\right),\quad \bar{v} = \max \left(1,\frac{1}{a_1}\right),
\end{align*}
then the hypothesis \eqref{H1} and \eqref{H2} are satisfied. 
If $d_1=d_2=1$ and $\alpha_1=\alpha_2=1$, then by \thref{prop: lower bed} (or by \thref{rem: lower bound equal diffusion}), 
$$
(u+k_1)(v+k_2)\geq \min\left( (\underline u+k_1)k_2, k_1(\underline v+k_2)\right),
\quad \forall k_1, k_2>0.
$$
Recall the maximum principle in Theorem 1.1 in \cite{JDE-16}:
$$
\min\left(\frac{k_2}{a_2}, \frac{k_1}{a_1}\right)\leq k_2 u+k_1 v\leq \max\left(k_1, k_2\right),
$$
then we have 
$$
uv+\max\left( k_1, k_2\right)\geq uv+k_2 u+k_1 v\geq \min\left( k_2\underline u, k_1\underline v\right).
$$
Under the bistable condition $a_1,a_2>1$, we derive the following ``trivial'' lower bound by taking $k_1=k_2$,
$$
uv+k_1\geq k_1\min(\frac{1}{a_1}, \frac{1}{a_2})
\Rightarrow uv\geq 0.
$$
According to \eqref{eqn: upper bound of p general case}, letting $\alpha_1=\alpha_2=m_1=m_2=1$ leads to
\begin{equation}\label{eqn: uv upper bound eg 1}
\sqrt{u(x)v(x)}\le
\frac{1}{2}
\max \left( \bar{u}, \bar{v}\right)
\frac{\max \left(d_1,d_2\right)}{\min \left(d_1,d_2\right)}, \quad x\in\mathbb{R}
\end{equation}
or 
\begin{equation}\label{eqn: uv upper bound eg 2}
u(x)v(x)\le
\frac{1}{4}
\left(\max \left( \bar{u}, \bar{v}\right)\right)^2
\left(\frac{\max \left(d_1,d_2\right)}{\min \left(d_1,d_2\right)}\right)^2, \quad x\in\mathbb{R}.
\end{equation}
For the equal diffusion case $d_1=d_2=1$ with the bistable condition $a_1,a_2>1$, \eqref{eqn: uv upper bound eg 2} is simplified to
\begin{equation}\label{eqn: uv upper bound eg 3}
u(x)v(x)\le\frac{1}{4}, \quad x\in\mathbb{R}.
\end{equation}
If we further consider the boundary conditions in the $(\textbf{e}_2,\textbf{e}_4)$-waves (also $(\textbf{e}_3,\textbf{e}_4)$-waves) or the $(\textbf{e}_4,\textbf{e}_4)$-waves, the upper bound  given by \eqref{eqn: uv upper bound eg 3} is optimal for the case $a:=a_1=a_2>1$ since as $a\to1^{+}$, we have
\begin{equation}
(u^{\ast},v^{\ast})=\left(\displaystyle\frac{1-a_1}{1-a_1\,a_2},\displaystyle\frac{1-a_2}{1-a_1\,a_2}\right)=\left(\displaystyle\frac{1}{1+a},\displaystyle\frac{1}{1+a}\right)\to\left(1/2,1/2\right).
\end{equation}

\section{Proofs of \thref{prop: lower bed} and \thref{prop: upper bed}}\label{Sec: Proofs}


\begin{proof}[Proof of \thref{prop: lower bed}]
We first rewrite the inequality $d_i(u_i)''+\theta(u_i)'+u_i^{l_i}f_i\le0$ in \eqref{eqn:upper}. If $u(x)\geq 0$, then for any $k>0$, a straightforward calculation gives
\begin{align*}  
(\ln (u(x)+k))'&=\frac{u'(x)}{u(x)+k},\\ 
(\ln (u(x)+k))''&=\frac{u''(x)}{u(x)+k}-\frac{(u'(x))^2}{(u(x)+k)^2}.
\end{align*}
Hence we divide the inequality by  \bluee{$u_i+k_i>0$} with $k_i>0$ to arrive at
\begin{align*}
\notag
 d_i(\ln (u_i+k_i))''+d_i\frac{((u_i)')^2}{(u_i+k_i)^2}+\theta\,(\ln (u_i+k_i))'+\frac{u_i^{l_i}}{u_i+k_i}f_i
\leq 0.
\end{align*}
Thus $(U_i)_{i=1}^n:=(\ln(u_i+k_i))_{i=1}^n$ satisfies the following inequalities:
\begin{align} \label{eqn: differential ineq lower bound ln(u+k)}
d_i U_i''+\theta\,U_i'+\frac{u_i^{l_i}}{u_i+k_i}f_i\leq 0,\quad i=1,\cdots,n.
\end{align}
For any   $(\alpha_i)_{i=1}^n\in (\mathbb{R}^+)^n$, let
\begin{equation*}
p(x) =\sum_{i=1}^{n} \alpha_i U_i,
\quad 
q(x) =\sum_{i=1}^{n} \alpha_i\,d_i U_i,
\end{equation*}
then the above inequality \eqref{eqn: differential ineq lower bound ln(u+k)} reads as
\begin{equation}\label{eq:pq}
q''+\theta p'+F\le0, 
\quad F:=\displaystyle\sum_{i=1}^{n} \frac{\alpha_i\,u_i^{l_i}}{u_i+k_i}\,f_i(u_1,\cdots,u_n).
\end{equation}

We are going to derive a lower bound for 
$$
q=\sum_{i=1}^n\alpha_i d_i U_i=\sum_{i=1}^{n} \alpha_i\,d_i\ln(u_i(x)+k_i),
$$
and hence a lower bound for $\prod_{i=1}^n(u_i+k_i)^{d_i\alpha_i}$.
The idea is similar as in the papers  \cite{NBMP-Discrete,JDE-16,CPAA-16,DCDS-B-18,NBMP-n-species,DCDS-A-17}, namely we are going to determine  three parameters
$$
\lambda_1, \quad \eta,\quad\lambda_2
$$
to construct an N-barrier consisting of three hypersurfaces 
\begin{align*}
Q_1:=\{(u_i)_{i=1}^n\,|\, q=\lambda_1\},
\quad P:=\{(u_i)_{i=1}^n\,|\, p=\eta\},
\quad Q_2:=\{(u_i)_{i=1}^n\,|\, q=\lambda_2\},
\end{align*}
such that the following inclusion relations hold:
\begin{align*}\label{eqn: Q1<P<Q2<R uv} 
&\mathcal{Q}_{1}:=\{(u_i)_{i=1}^n\in([0,\infty))^n
\,|\, q\leq \lambda_1\}
 \subset \mathcal{P}:=\{(u_i)_{i=1}^n\in([0,\infty))^n
 \,|\, p\leq \eta\}
\\
& \subset \mathcal{Q}_{2}:=\{(u_i)_{i=1}^n\in([0,\infty))^n
\,|\, q\leq\lambda_2\}
  \subset  \underaccent\bar{\mathcal{R}}
= \big\{ (u_i)_{i=1}^n\in([0,\infty))^n \; |\; \sum_{i=1}^{n}\frac{\displaystyle u_i}{\displaystyle\underaccent\bar{u}_i}\le 1\big\}.
\end{align*}
It will turn out that if  $\lambda_1$, $\eta$, and $\lambda_2$ are given respectively by \eqref{eqn: lambda1 lower bound}, \eqref{eqn: eta lower bound}, and \eqref{eqn: lambda2 lower bound}, then $\lambda_1$ determines a lower bound of $q(x)$: $q(x)\geq\lambda_1$, which is exactly \eqref{ineq:lower}.

More precisely, we follow the steps as in  \cite{NBMP-Discrete,JDE-16,CPAA-16,DCDS-B-18,NBMP-n-species,DCDS-A-17} to  determine $\lambda_2$, $\eta$, $\lambda_1$ such that the above inclusion relations $\mathcal{Q}_1\subset\mathcal{P}\subset\mathcal{Q}_2\subset\underaccent\bar{\mathcal{R}}$ hold:
\begin{enumerate}[$(i)$]
  \item \textbf{\underline{Determine $\lambda_2$}}  
  The hypersurface $Q_2$ intersects the  $u_j$-axis: $\{(u_i)_{i=1}^n\,|\,u_i=0,\,\forall i\neq j\}$ at the point
    \begin{equation*}\label{eqn: intercept lambda2}
\Bigl\{(u_i)_{i=1}^n\,|\, u_i=0,\,\forall i\neq j,
\quad   u_{2,j}=e^{ \frac{\lambda_2-\sum_{i=1,i\neq j}^{n} \alpha_i\,d_i\ln k_i}{\alpha_j\,d_j}}-k_j\Bigr\}.
   \end{equation*}
If $u_{2,j}\leq \underaccent\bar{u}_j$, $\forall j=1,\cdots,n$, then by the monotonicity of the function $\ln(\cdot+k)$, $\mathcal{Q}_2\subset \underaccent\bar{\mathcal{R}}$.   
That is, $\mathcal{Q}_2\subset \underaccent\bar{\mathcal{R}}$ if $\lambda_2$ is chosen as in \eqref{eqn: lambda2 lower bound}:
    \begin{equation}\nonumber
  \lambda_2=\min_{1\le j \le n} \Bigl(\alpha_j d_j \ln (\underaccent\bar{u}_j+k_j)+\sum_{i=1, i\neq j}^{n}\alpha_i d_i\ln k_i\Bigr).
   \end{equation} 

\item \textbf{\underline{Determine $\eta$}}  
As above, the hypersurface $P$ intersects the $u_j$-axis  at
    \begin{equation*}
\Bigl\{(u_i)_{i=1}^n\,|\, u_i=0,\,\forall i\neq j,
\quad   u_{0,j}=e^{ \frac{\eta-\sum_{i=1,i\neq j}^{n} \alpha_i\, \ln k_i}{\alpha_j}}-k_j\Bigr\}.
   \end{equation*}
If $u_{0,j}\leq u_{2,j}$, $\forall j=1,\cdots,n$, then $\mathcal{P}\subset \mathcal{Q}_2$ and the hypersurface $Q_2$ is above the hypersurface $P$. That is, $\mathcal{P}\subset \mathcal{Q}_2$ if $\eta$ is chosen as  in \eqref{eqn: eta lower bound}:
    \begin{equation}\nonumber
  \eta=\min_{1\le j \le n} \frac{1}{d_j}\Bigl(\lambda_2-\sum_{i=1, i\neq j}^{n}\alpha_i(d_i-d_j)\ln k_i\Bigr).
   \end{equation}
  
  \item \textbf{\underline{Determine $\lambda_1$}} \bluee{Replacing $\lambda_2$ by $\lambda_1$ in step $(i)$},  the $u_j$-intercept of the hypersurface $Q_1$ is given by 
$u_{1,j}=e^{\frac{\lambda_1-\sum_{i=1,i\neq j}^{n} \alpha_i\,d_i\ln k_i}{\alpha_j\,d_j}}-k_j$.
Hence if we take $\lambda_1$ as  in \eqref{eqn: lambda1 lower bound}:
    \begin{equation}\nonumber
 \lambda_1=\min_{1\le j \le n} \Bigl(\eta\,d_j+\sum_{i=1, i\neq j}^{n}\alpha_i(d_i-d_j)\ln k_i\Bigr),
   \end{equation}  
   then $u_{1,j}\leq u_{0,j}$, $\forall j=1,\cdots,n$ and hence $\mathcal{Q}_1\subset \mathcal{P}$. 
   \end{enumerate}
   
   We now show $q(x)\ge \lambda_1$, $x\in\mathbb{R}$ by a contradiction argument. 
   Suppose by contradiction that there exists $z\in\mathbb{R}$ such that $q(z)<\lambda_1$. Since $u_i(x)\in C^2(\mathbb{R})$ $(i=1,\cdots,n)$ and  $(u_1,u_2,\cdots,u_n)(\pm\infty)=\textbf{e}_{\pm}$, we may assume $\displaystyle\min_{x\in\mathbb{R}} q(x)=q(z)$. 
We denote respectively by $z_2$ and $z_1$ the first points at which the solution trajectory 
$\{(u_i(x))_{i=1}^n\,|\,x\in\mathbb{R}\}$ intersects the hypersurface $Q_2$ when $x$ moves from $z$ towards $\infty$ and $-\infty$. 
For the case where $\theta\leq0$, we integrate \eqref{eq:pq} with respect to $x$ from $z_1$ to $z$ and obtain
\begin{equation}\label{eqn: integrating eqn}
q'(z)-q'(z_1)+\theta\,(p(z)-p(z_1))+\int_{z_1}^{z}F(u_1(x),\cdots,u_n(x))\,dx\leq0.
\end{equation}
We also have the following facts from the construction of the hypersurfaces $Q_1, Q_2, P$:
\begin{itemize}
  \item $q'(z)=0$ because of $\displaystyle\min_{x\in\mathbb{R}} q(x)=q(z)$;
  \item $q(z_1)=\lambda_2$ because of $(u_i(z_1))_{i=1}^n\in Q_2$. 
  \item $q'(z_1)< 0$ because   $z_1$ is the first point
  for $q(x)$ taking the value $\lambda_2$ when $x$ moves from $z$ to $-\infty$, such that $q(z_1+\delta)< \lambda_2$ for $z-z_1>\delta>0$;
  \item $p(z)<\eta$ since $(u_i(z))_{i=1}^n$ is below the hypersurface $P$;
  \item $p(z_1)>\eta$ since $(u_i(z_1))_{i=1}^n$ is above the hypersurface $P$;
  \item $F(u_1(x),\cdots,u_n(x))=\displaystyle\sum_{i=1}^{n} \frac{\alpha_i\,u_i^{l_i}}{u_i+k_i}\,f_i(u_1,\cdots,u_n)\geq0$, $\forall x\in [z_1,z]$.
  Indeed, since $(u_i(z_1))_{i=1}^n\in Q_2\subset\mathcal{Q}_2\subset\underaccent\bar{\mathcal{R}}$ and $(u_i(z))_{i=1}^n\in \mathcal{Q}_1\subset\underaccent\bar{\mathcal{R}}$, we derive that $F(u_1(x),\cdots,u_n(x))|_{x\in[z_1,z]}\geq0$ by the hypothesis \eqref{H1}. 
\end{itemize}
We hence have the following inequality from the above facts when $\theta\leq 0$
\begin{equation*}
q'(z)-q'(z_1)+\theta\,(p(z)-p(z_1))+\int_{z_1}^{z}F(u_1(x),\cdots,u_n(x))\,dx>0,
\end{equation*}
which contradicts \eqref{eqn: integrating eqn}. Therefore when $\theta\leq0$, $q(x)\geq \lambda_1$ for $x\in \mathbb{R}$. For the case where $\theta\geq0$, we simply integrate \eqref{eq:pq}  with respect to $x$ from $z$ to $z_2$ to arrive at
\begin{equation*}\label{eqn: eqn by integrate from z to z2}
q'(z_2)-q'(z)+\theta\,(p(z_2)-p(z))+\int_{z}^{z_2}F(u_1(x),\cdots,u_n(x))\,dx\leq0.
\end{equation*}
Then we apply the facts that $q'(z_2)> 0$, $q'(z)=0$, $p(z_2)>\eta$, $p(z)<\eta$ and $F(u_1(x),\cdots,u_n(x))|_{x\in[z,z_2]}\geq 0$, as well as  a similar contradiction argument as above, to derive $q\geq\lambda_1$.



\end{proof}

\begin{proof}[Proof of \thref{prop: upper bed}] 
We prove \thref{prop: upper bed} in a similar manner to the proof of \thref{prop: lower bed}. 
We first rewrite the inequality $d_i(u_i)''+\theta(u_i)'+u_i^{l_i}f_i\geq 0$ in \eqref{eqn:lower}.
A straightforward calculation shows
\begin{align*}   
&(u^m )'=m\,u^{m-1} u',\\[1ex] 
&(u^m)''=m\left( (m-1)\,u^{m-2}(u')^2+u^{m-1}u''(x)\right).
\end{align*}
Hence we multiply the inequality by $m_i\,u^{m_i-1}(x)$ to arrive at
\begin{align*}
&d_i ({u_i}^{m_i})'' -d_i\,m_i(m_i-1)\,{u_i}^{m_i-2}({u_i}')^2+\theta\,( {u_i}^{m_i})'+m_i\,{u_i}^{m_i-1}u_i^{l_i}\,f_i 
 \geq 0.
\end{align*} 
For notational simplicity, we will adopt the same notations as in  the proof of \thref{prop: lower bed}.
Since $u_i\geq 0$, $\forall i=1,\cdots,n$, for any $(m_i)_{i=1}^n\in ([1,\infty))^n$, the vector field $(U_i)_{i=1}^n:=(u_i^{m_i})_{i=1}^n$ satisfies the following inequalities
\begin{align}  \label{eqn: differential ineq upper bound u^m}
 d_i U_i''+\theta\,U_i'+m_i\,{u_i}^{m_i-1}u_i^{l_i}\,f_i\geq 0, \quad \forall i=1,\cdots,n.
\end{align}
For any   $(\alpha_i)_{i=1}^n\in (\mathbb{R}^+)^n$, 
$p(x)=\sum_{i=1}^{n} \alpha_i\,U_i$ and
$q(x)=\sum_{i=1}^{n} \alpha_i\,d_i U_i$
satisfy
\begin{equation}\label{eqn: ODE for p and q upper bound}
q''+\theta\,p'+F\ge0,\quad  F:=\displaystyle\sum_{i=1}^{n} \alpha_i\,m_i\,{u_i}^{m_i-1}u_i^{l_i}\,f_i(u_1,u_2,\cdots,u_n).
\end{equation}

We are going to show the upper bound $q\leq \lambda_1$ by employing the N-barrier method as in the proof of Proposition~\ref{prop: lower bed}. 
That is, we are going to construct the three hyperellipsoids 
\begin{align*}
Q_1:=\{(u_i)_{i=1}^n\,|\, q=\lambda_1\},
\quad P:=\{(u_i)_{i=1}^n\,|\, p=\eta\},
\quad Q_2:=\{(u_i)_{i=1}^n\,|\, q=\lambda_2\},
\end{align*}
such that the following inclusion relations hold:
\begin{align*}\label{eqn: Q1<P<Q2<R uv} 
&\mathcal{Q}_{1}:=\{(u_i)_{i=1}^n\in([0,\infty))^n
\,|\, q\geq \lambda_1\}
 \supset \mathcal{P}:=\{(u_i)_{i=1}^n\in([0,\infty))^n
 \,|\, p\geq \eta\}
\\
& \supset \mathcal{Q}_{2}:=\{(u_i)_{i=1}^n\in([0,\infty))^n
\,|\, q\geq\lambda_2\}
  \supset  \bar{\mathcal{R}}
= \big\{ (u_i)_{i=1}^n\in([0,\infty))^n \; |\; \sum_{i=1}^{n}\frac{\displaystyle u_i}{\displaystyle \bar{u}_i}\ge 1\big\},
\end{align*} 
and the upper bound $q\leq\lambda_1$ follows by a contradiction argument.
More precisely, we take 
    \begin{equation}\label{eqn: lambda2 upper bound uv}
  \lambda_2=\max_{1\le i\le n} \alpha_i\,d_i(\bar{u}_i)^{m_i},
   \end{equation}
   such that the $u_j$-intercept of the hyperellipsoid $Q_2$ 
  \begin{equation*}
    u_{2,j}=\left(\frac{\lambda_2}{\alpha_j\,d_j}\right)^{1/m_i}\geq \bar u_j, \quad j=1,2,\cdots,n.
   \end{equation*}
Then we take 
    \begin{equation}\label{eqn: eta upper bound uv}
  \eta= \frac{\lambda_2}{\displaystyle\min_{1\le i \le n} d_i},
   \end{equation}
   such that the $u_j$-intercept of the hyperellipsoid $P$ 
  \begin{equation*}\label{eqn: intercept eta upper bound}
  u_{0,j}=\left(\frac{\eta}{\alpha_j}\right)^{1/m_i}\geq u_{2,j}, \quad j=1,2,\cdots,n.
   \end{equation*}
   Finally we take 
    \begin{equation}\label{eqn: lambda1 upper bound uv}
  \lambda_1=\eta\max_{1\le i \le n} d_i
   \end{equation}
   such that  the $u_j$-intercept of the hyperellipsoid $Q_1$ 
  \begin{equation*}\label{eqn: intercept lambda1 upper bound}
  u_{1,j}=\left(\frac{\lambda_1}{\alpha_j\,d_j}\right)^{1/m_i}\geq u_{0,j}, \quad j=1,2,\cdots,n.
   \end{equation*} 
Combining \eqref{eqn: lambda2 upper bound uv}, \eqref{eqn: eta upper bound uv}, and \eqref{eqn: lambda1 upper bound uv}, we have 
\begin{equation}
\lambda_1=\left(\max_{1\le i \le n} \alpha_i\,d_i(\bar{u}_i)^{m_i}\right) 
\frac{\displaystyle\max_{1\le i \le n} d_i}{\displaystyle\min_{1\le i \le n} d_i}.
\end{equation} 
We follow exactly the same contradiction argument to prove $q(x)\le\lambda_1$ for $x\in\mathbb{R}$ 
as in the proof of \thref{prop: lower bed}, which is omitted here.
 Since $(\alpha_i)_{i=1}^n\in (\mathbb{R}^+)^n$  is arbitrary, $q(x)=\displaystyle\sum_{i=1}^{n} \alpha_i\,d_i(u_i(x))^{m_i}\le\lambda_1$  implies the upper bound   \eqref{bound:upper}.
Now we use the inequality of arithmetic and geometric means 
to obtain
\begin{align}
\label{eqn: imply of AM?GM inequality} 
\sum_{i=1}^{n} \alpha_i(u_i(x))^{m_i}
\ge n\Bigl(\prod_{i=1}^{n} \alpha_i(u_i(x))^{m_i}\Bigr)^{\frac{1}{n}} 
\ge  n\Bigl(\prod_{i=1}^{n}\alpha_i\Bigr)^{\frac{1}{n}}  \prod_{i=1}^{n}(u_i(x))^{\frac{m_i}{n}},
\end{align}
which together with \eqref{bound:upper} yields  \eqref{eqn: upper bound of p general case}.

\end{proof}

\subsubsection*{Acknowledgements}
The authors are grateful to the anonymous referees for many helpful comments and valuable suggestions on this paper. L.-C. Hung thanks for the hospitality he received from KIT while visiting KIT. The research of L.-C. Hung is partly supported by the grant 106-2115-M-011-001-MY2 of Ministry of Science and Technology, Taiwan.
 





\begin{thebibliography}{9999}  
 
\bibitem{Ahmad&Lazer91}
{\sc S.~Ahmad and A.~C. Lazer}, {\em An elementary approach to traveling front
  solutions to a system of {$N$} competition-diffusion equations}, Nonlinear
  Anal., 16 (1991), pp.~893--901.

\bibitem{NBMP-Discrete}
{\sc C.-C. Chen, T.-Y. Hsiao, and L.-C. Hung}, {\em Discrete n-barrier maximum
  principle for a lattice dynamical system arising in competition models}, to
  appear in Discrete Contin. Dyn. Syst. A.

\bibitem{JDE-16}
{\sc C.-C. Chen and L.-C. Hung}, {\em A maximum principle for diffusive
  lotka-volterra systems of two competing species}, J. Differential Equations,
  261 (2016), pp.~4573--4592.

\bibitem{CPAA-16}
\leavevmode\vrule height 2pt depth -1.6pt width 23pt, {\em Nonexistence of
  traveling wave solutions, exact and semi-exact traveling wave solutions for
  diffusive {L}otka-{V}olterra systems of three competing species}, Commun.
  Pure Appl. Anal., 15 (2016), pp.~1451--1469.

\bibitem{DCDS-B-18}
\leavevmode\vrule height 2pt depth -1.6pt width 23pt, {\em An n-barrier maximum
  principle for elliptic systems arising from the study of traveling waves in
  reaction-diffusion systems}, Discrete Contin. Dyn. Syst. B, 22 (2017),
  pp.~1--19.

\bibitem{NBMP-n-species}
{\sc C.-C. Chen, L.-C. Hung, and C.-C. Lai}, {\em An n-barrier maximum
  principle for autonomous systems of n species and its application to problems
  arising from population dynamics}, Commun. Pure Appl. Anal., 18 (2019),
  pp.~33--50.

\bibitem{DCDS-A-17}
{\sc C.-C. Chen, L.-C. Hung, and H.-F. Liu}, {\em N-barrier maximum principle
  for degenerate elliptic systems and its application}, Discrete Contin. Dyn.
  Syst. A, 38 (2018), pp.~791--821.

\bibitem{Fei&Carr03}
{\sc N.~Fei and J.~Carr}, {\em Existence of travelling waves with their minimal
  speed for a diffusing {L}otka-{V}olterra system}, Nonlinear Anal. Real World
  Appl., 4 (2003), pp.~503--524.

\bibitem{Hou&Leung08}
{\sc X.~Hou and A.~W. Leung}, {\em Traveling wave solutions for a competitive
  reaction-diffusion system and their asymptotics}, Nonlinear Anal. Real World
  Appl., 9 (2008), pp.~2196--2213.

\bibitem{Kan-on95}
{\sc Y.~Kan-on}, {\em Parameter dependence of propagation speed of travelling
  waves for competition-diffusion equations}, SIAM J. Math. Anal., 26 (1995),
  pp.~340--363.

\bibitem{Kan-on97Fisher-Monostable}
\leavevmode\vrule height 2pt depth -1.6pt width 23pt, {\em Fisher wave fronts
  for the {L}otka-{V}olterra competition model with diffusion}, Nonlinear
  Anal., 28 (1997), pp.~145--164.

\bibitem{Kanel06}
{\sc J.~I. Kanel}, {\em On the wave front solution of a competition-diffusion
  system in population dynamics}, Nonlinear Anal., 65 (2006), pp.~301--320.

\bibitem{Kanel&Zhou96}
{\sc J.~I. Kanel and L.~Zhou}, {\em Existence of wave front solutions and
  estimates of wave speed for a competition-diffusion system}, Nonlinear Anal.,
  27 (1996), pp.~579--587.

\bibitem{Leung&Feng09}
{\sc A.~W. Leung, X.~Hou, and W.~Feng}, {\em Traveling wave solutions for
  lotka-volterra system re-visited}, Discrete \& Continuous Dynamical
  Systems-B, 15 (2011), pp.~171--196.

\bibitem{Leung08}
{\sc A.~W. Leung, X.~Hou, and Y.~Li}, {\em Exclusive traveling waves for
  competitive reaction-diffusion systems and their stabilities}, J. Math. Anal.
  Appl., 338 (2008), pp.~902--924.

\bibitem{Murray93Mbiology}
{\sc J.~D. Murray}, {\em Mathematical biology}, vol.~19 of Biomathematics,
  Springer-Verlag, Berlin, second~ed., 1993.

\bibitem{tang1980pfc}
{\sc M.~Tang and P.~Fife}, {\em {Propagating fronts for competing species
  equations with diffusion}}, Archive for Rational Mechanics and Analysis, 73
  (1980), pp.~69--77.

\bibitem{Volpert94}
{\sc A.~I. Volpert, V.~A. Volpert, and V.~A. Volpert}, {\em Traveling wave
  solutions of parabolic systems}, vol.~140 of Translations of Mathematical
  Monographs, American Mathematical Society, Providence, RI, 1994.
\newblock Translated from the Russian manuscript by James F. Heyda.

              
\end{thebibliography}



\end{document}